\documentclass{article}
\usepackage{amssymb}

\def\ac{a^{+}}
\def\al{\alpha}
\def\be{\beta}
\def\la{\lambda}
\def\ga{\gamma}
\def\si{\sigma}
\def\bal{\bar{\alpha}}
\def\bbe{\bar{\beta}}
\def\nn{\nonumber \\}
\def\ha{\frac{1}{2}}
\def\vp{\varphi}
\def\ve{\varepsilon}

\begin{document}

\begin{titlepage}

\vspace{4em}
\begin{center}

{\Large{\bf  The Eigenfunctions of the 
    \boldmath $q$-Harmonic Oscillator on the Quantum Line}}

\vskip 5em

{{\bf
H. Grosse, S. Schraml}}

\vskip 1em

   Insititut f\"ur Theoretische Physik\\
   der Universit\"at Wien\\
   Boltzmanngasse 5, A-1090 Wien, AUSTRIA\\[1em]

\end{center}

\vspace{2em}

\begin{abstract}
We construct a complete set of eigenfunctions of the 
$q$-deformed harmonic oscillator on the quantum line. 
In particular the eigenfunctions corresponding to the non-Fock 
part of the spectrum will be constructed.
\end{abstract}

\vfill
\noindent \hrule\vskip.1cm
\hbox{{\small{\it e-mail: }}{
\small\quad harald.grosse@univie.ac.at,\quad schraml@ap.univie.ac.at}}
\end{titlepage}\vskip.2cm

\newpage
\setcounter{page}{1}

\section{Introduction}

In this paper we will construct a complete set of eigenfunctions 
for the $q$-deformed harmonic oscillator on the $q$-deformed line. 
We define the $q$-deformed oscillator \cite{Macfarlane, Biedenharn} 
as the unital $*$-algebra generated by 
the element $a$ and its conjugate $\ac$, subject to the relation
\begin{equation}
  \label{qOszRel}
  a\ac - q^{-2}\ac a=1.
\end{equation}
We take $q>1$, since this will be the case for the realization of these 
operators, which we consider below.

The spectrum of the Hamilton operator $H\equiv\ac a$ consists of two parts 
\cite{ChGrPr}. 
A bounded part (Fock representation) and an unbounded part:
\begin{equation}
\label{qH-Spec}
  Spec\,(H)=\left\{
    \begin{array}{ll}
      \frac{1-q^{-2n}}{1-q^{-2}},& n\in\mathbb{N}\\
      \frac{1+q^{2\ga-2m}}{1-q^{-2}},&m\in\mathbb{Z},\; \ga\in\mathbb{R}
    \end{array}
  \right.
\end{equation}
Both parts of the spectrum have an accumulation point at $\frac{1}{1-q^{-2}}$. 
The Fock representation is a lowest weight representation with $\ac$
acting as raising operator and $a$ acting as lowering operator.
For the second, i.e. the unbounded, part of the spectrum $a$ is  
the raising operator.

It is known, that on the $q$-deformed line the $q$-deformed Hermite 
polynomials are related to the eigenfunctions 
corresponding to the Fock representation 
\cite{LRW,hinterding}. However, these functions are not complete in the respective 
Hilbert space of square integrable functions. One has to consider 
the eigenfunctions related to the unbounded part of the spectrum as well. 
This we will do by using results of Ciccoli et.al. \cite{Koorn}.

\section{Representation on the $q$-deformed line} 

We will consider a realization of the $q$-oscillator on the
$q$-deformed real line $\mathbb{R}_q$, which is defined as being generated 
by the operators $X,P,U$ with commutation relations
\begin{eqnarray}
  \label{XPURel}
  &q^{\ha}XP-q^{-\ha}PX=iU&\nn
  &UX=q^{-1}XU,\qquad UP=qPU&
\end{eqnarray}
and the following conjugation:
\begin{equation}
  X^+=X,\qquad P^+=P,\qquad U^+=U^{-1}
\end{equation}

A realization of the $q$-oscillator on $\mathbb{R}_q$ is 
given by \cite{LRW}:
\begin{eqnarray}
  \label{qosz-realization}
  a&=&\al U^{-2} + \be U^{-1}P\nn
  \ac&=&\bal U^2 + \bbe PU
\end{eqnarray}
With $\al,\be\in\mathbb{C}$, such that
\begin{equation}
  \al\bal=\frac{1}{1-q^{-2}}=\frac{q}{\la}\quad\mbox{and}\quad \al\bbe=\bal\be.
\end{equation}
The second relation implies
$\frac{\al}{\be}=\frac{\bal}{\bbe}\in\mathbb{R}$; we define
\begin{equation}
  \label{ga-def}
  q^{-\ga}\equiv\frac{\al}{\be}=\frac{\bal}{\bbe}
\end{equation}

The algebra (\ref{XPURel}) of the $q$-deformed real line
$\mathbb{R}_q$ can be realized by operators acting on functions of
one variable:
\begin{equation}
  \label{XPU-realization}
  Xf(x)=xf(x);\qquad Pf(x)=-iD_qf(x),\qquad Uf(x)=q^{-\ha}f(q^{-1}x);
\end{equation}
where a $q$-derivative in the following form has been used:
\begin{equation}
  \label{Dq-def}
  D_qf(x)\equiv\frac{f(qx)-f(q^{-1}x)}{x(q-q^{-1})}.
\end{equation}
The algebra acts on functions on a 'lattice' $\xi q^{n}$, 
$\xi\in\mathbb{R}, n\in\mathbb{Z}$.

The scalar product can be defined in terms of the Jackson
integral:
\begin{equation}
  (f,g)\sim\sum_{n\in\mathbb{Z}} \overline{f(q^n)}g(q^n)q^n
\end{equation}

We will use the following notation:
\begin{equation}
  (a;q)_n\equiv \prod_{i=0}^{n-1} (1-aq^i);\qquad (a;q)_{\infty}\equiv \lim_{n\to\infty}(a;q)_n
\end{equation}
and define the $q$-exponential function:
\begin{eqnarray}
  e_q(x)&\equiv&\frac{1}{(x;q^{-2})_{\infty}}\nn
  D_qe_q(cx)&=&c\frac{q}{\la}e_q(qcx),\qquad c\in\mathbb{C}
\end{eqnarray}

Using the relations (\ref{qosz-realization}) and
(\ref{XPU-realization}), it is easily seen, that the ground state
of the Fock representation, i.e. the state satisfying
$a\psi_0(x)=0$, is given by:
\begin{equation}
  \psi_0(x)\equiv Ne_q(-i\frac{\al}{\be}\la q^{-\ha}x),
\end{equation}
where $N$ is a normalization constant.

The Hamilton operator $H\equiv \ac a$ is in terms of $D_q$ and
$U$:
\begin{equation}
  H=\ac a=\al\bal-\be\bbe D_q^2-i\al\bbe(U+qU^{-1})D_q
\end{equation}
With this, and the action (\ref{XPU-realization}), (\ref{Dq-def}) 
 of $D_q$ and $U$ on functions, the equation
\begin{equation}
  H f(x)=Ef(x)
\end{equation}
for the eigenfunctions becomes a difference equation:
\begin{eqnarray}
\label{f-EV-eqn}
  Ex^2\la^2f(x)&=&f(x)\left\{\al\bal x^2\la^2+\be\bbe(q+q^{-1})\right\}\nn
  &&+f(q^2x)\left\{-q^{-1}\be\bbe-i\al\bbe q^{\ha}x\la\right\}\nn
  &&+f(q^{-2}x)\left\{-q\be\bbe+i\al\bbe q^{\ha}x\la\right\}
\end{eqnarray}

With the definition
\begin{equation}
  E=\frac{1+\ve}{1-q^{-2}}
\end{equation}
the equation (\ref{f-EV-eqn}) for the eigenstates $f(x)=\psi_0(x)g(x)$
with eigenvalue $E$ becomes:
\begin{eqnarray}
\label{EV-eq}
  0&=&g(x)\left\{q+q^{-1}-\ve q^{-2\ga}\la^2x^2\right\}\nn
  &&-q^{-1}g(q^2x)-qg(q^{-2}x)\left\{1+q^{-2\ga-1}\la^2x^2\right\}
\end{eqnarray}
Where we now use $\gamma$, Eqn. (\ref{ga-def}), instead of $\al$ and $\beta$.

\section{Orthonormal basis}

To solve this equation, we use the basic hypergeometric series
${_1}\vp_1$ \cite{GaRa}. The function $f(z)={_1}\vp_1(a;c;q,z)$ satisfies:
\begin{equation}
\label{phi11-recursion}
  (c-az)f(qz)+(-(c+q)+z)f(z)+qf(z/q)=0
\end{equation}
We define:
\begin{eqnarray}
  \vp_e(x)&\equiv& {_1}\vp_1(-\ve^{-1};q^{-2};q^{-4},\ve q^{-2\ga-3}\la^2x^2)\nn
  \vp_o(x)&\equiv& {_1}\vp_1(-q^{-2}\ve^{-1};q^{-6};q^{-4},\ve q^{-2\ga-5}\la^2x^2).
\end{eqnarray}
Due to the relation (\ref{phi11-recursion}) we find that the following functions
solve the Eqn. (\ref{EV-eq}):
\begin{equation}
\label{solutions}
  g(x)=\vp_e(x),\qquad g(x)=x\vp_o(x)
\end{equation}
These two solutions correspond to parts of the spectrum (\ref{qH-Spec}), which are 
numbered by even and odd numbers respectively.

It is possible, to combine the two solutions (\ref{solutions}) to a function, 
that yields the whole spectrum \cite{Koorn}:
\begin{eqnarray}
\label{gesLsg}
  \lefteqn{(-\ve)^k{_2}\vp_1\left(-\frac{1}{\ve q^2},-\frac{1}{\ve};0;q^{-4},-\frac{q^{4(k-1)}}{c}\right)}\\
  &&=C_e{_1}\vp_1\left(-\frac{1}{\ve};q^{-2};q^{-4},\frac{\ve c}{q^{4k+2}}\right)+C_oq^{-2k}{_1}\vp_1\left(-\frac{1}{\ve q^2};q^{-6};q^{-4},\frac{\ve c}{q^{4k+4}}\right)\nonumber
\end{eqnarray}
with
\begin{eqnarray}
  C_e&=&\frac{(-\ve^{-1}q^{-2},\ve^{-1}c^{-1}q^{-4},\ve c;q^{-4})_{\infty}}{(q^{-2},-c,-c^{-1}q^{-4};q^{-4})_{\infty}}\nn
  C_o&=&\frac{(-\ve^{-1},\ve^{-1}c^{-1}q^{-6},\ve cq^2;q^{-4})_{\infty}}{(q^{-6},-c,-c^{-1}q^{-4};q^{-4})_{\infty}}
\end{eqnarray}

We will discuss some properties of these functions and consider  
later the relation to the lattice and Hilbert space coming from representations of 
the algebra Eq. (\ref{XPURel}). 
First we consider the Fock representation ($\ve=-q^{-2p},\, p\in\mathbb{N}$). 
For the lattice points $x=\xi q^{2n}$ one has, using the results of
\cite{Koorn} (Theorem 4.1) in the case $\ve=-q^{-4p}$:
\begin{equation}
  \vp_e^p(n)={_1}\vp_1(q^{4p};q^{-2};q^{-4},-q^{-4p} q^{-2\ga-3}\la^2\xi^2q^{4n})
\end{equation}
with $c\equiv q^{-2\ga-1}\xi^2\la^2$ one obtains:
\begin{equation}
\label{FN1}
  \sum_{n=-\infty}^{+\infty}\frac{\vp_e^r(n)\vp_e^s(n)}{(-cq^{4n};q^{-4})_{\infty}}q^{2n}=\delta_{rs}q^{4r}\frac{(q^{-4};q^{-4})_r(q^{-4},-cq^{-2},-q^{-2}c^{-1};q^{-4})_{\infty}}{(q^{-2};q^{-4})_r(q^{-2},-c,-q^{-4}c^{-1};q^{-4})_{\infty}}
\end{equation}
Notice, that a $q$-exponential function turns up as measure 
under the Jackson integral, since $(x;q)_{\infty}(-x;q)_{\infty}=(x^2;q^2)_{\infty}$. 
This is similar to the undeformed case, where $e^{-\frac{1}{2}x^2}e^{-\frac{1}{2}x^2}=e^{-x^2}$ 
leads to the orthogonality measure for the Hermite polynomials.

For $\ve=-q^{-4p-2}$ one obtains
\begin{equation}
  \vp_o^p(n)=q^{2n}{_1}\vp_1(q^{4p};q^{-6};q^{-4},-q^{-4p-2} q^{-2\ga-5}\la^2\xi^2q^{4n})
\end{equation}
and
\begin{equation}
\label{FN2}
  \sum_{n=-\infty}^{+\infty}\frac{\vp_o^r(n)\vp_o^s(n)}{(-cq^{4n};q^{-4})_{\infty}}q^{2n}=\delta_{rs}q^{4r}\frac{(q^{-4};q^{-4})_r(q^{-4},-cq^{-6},-q^{2}c^{-1};q^{-4})_{\infty}}{(q^{-6};q^{-4})_r(q^{-6},-c,-q^{-4}c^{-1};q^{-4})_{\infty}}
\end{equation}

Now we use (\ref{gesLsg}) to combine these two parts. 
For $\ve=-q^{-4n}$ one finds:
\begin{equation}
  C_e=(-c)^{-n}q^{4n^2-2n}(q^{-2};q^{-4})_n, \qquad C_o=0
\end{equation}
Furthermore Eqn. (\ref{gesLsg}) becomes
\begin{eqnarray}
\label{FL1}
  &&\left(\pm\sqrt{c}q^{-2k}\right)^{2n}{_2}\vp_1\left(q^{4n-2},q^{4n};0;q^{-4},-\frac{q^{-4}}{(\pm\sqrt{c}q^{-2k})^2}\right)\\
  &=&(-)^nq^{4n^2-2n}(q^{-2};q^{-4})_n{_1}\vp_1\left(q^{4n};q^{-2};q^{-4},-q^{-4n-2}(\pm\sqrt{c}q^{-2k})^2\right)\nonumber
\end{eqnarray}
For $\ve=-q^{-4n-2}$ we find
\begin{equation}
  C_e=0,\qquad C_o=(-c)^{-n}q^{4n^2+2n}(q^{-6};q^{-4})_n
\end{equation}
and
\begin{eqnarray}
\label{FL2}
  &&\left(\pm\sqrt{c}q^{-2k}\right)^{2n+1}{_2}\vp_1\left(q^{4n},q^{4n+2};0;q^{-4},-\frac{q^{-4}}{(\pm\sqrt{c}q^{-2k})^2}\right)\\
  &=&(-)^{n+1}\sqrt{c}q^{4n^2+2n-2k}(q^{-6};q^{-4})_n{_1}\vp_1\left(q^{4n};q^{-6};q^{-4},-q^{-4n-6}(\pm\sqrt{c}q^{-2k})^2\right)\nonumber
\end{eqnarray}
Using $m=2n$ in the first case $\ve=-q^{-4n}$ and $m=2n+1$ for $\ve=-q^{-4n-2}$
the left hand sides of the Eqns. (\ref{FL1}) and (\ref{FL2}) are identical:
\begin{equation}
\label{qHermiteDef}
  \tilde{h}_m(x)=x^m{_2}\vp_1\left(q^{2m-2},q^{2m};0;q^{-4},-\frac{q^{-4}}{x^2}\right)
\end{equation}
The function $\psi_0(x)\tilde{h}(x)$ is the eigenfunction corresponding 
to the eigenvalue $\ve=-q^{-2m}$. 
In both cases we obtain from the Eqns. (\ref{FN1}), (\ref{FN2})
\begin{equation}
  \sum_{k=-\infty}^{+\infty}\frac{\tilde{h}_m(\sqrt{c}q^{-2k})\tilde{h}_m(\sqrt{c}q^{-2k})}{(-cq^{-4k};q^{-4})_{\infty}}q^{-2k}=N_{c}\frac{(q^{-2};q^{-2})_m}{q^{-2m^2}}
\end{equation}
where
$N_{c}=\frac{(q^{-4},-cq^{-2},-c^{-1}q^{-2};q^{-4})_{\infty}}{(q^{-2},-c,-c^{-1}q^{-4};q^{-4})_{\infty}}$.
The functions $\tilde{h}_m$ for even $m$ are not orthogonal to functions $\tilde{h}_m$ 
with odd $m$. Since according to the definition (\ref{qHermiteDef})
the functions $\tilde{h}_m(x)$ are even and odd for even $m$ and
odd $m$ respectively, we extend the sum to negative $x$-values.
Then we have:
\begin{equation}
\label{Hermite-ortho}
  \sum_{k=-\infty,\si=\pm}^{+\infty}\frac{\tilde{h}_m(\si\sqrt{c}q^{-2k})\tilde{h}_n(\si\sqrt{c}q^{-2k})}{(-cq^{-4k};q^{-4})_{\infty}}q^{-2k}=2N_{c}\frac{(q^{-2};q^{-2})_m}{q^{-2m^2}}\delta_{mn},
\end{equation}
which is the well known orthogonality relation for the $q$-Hermite
II polynomials, that are known to be related to the $q$-oscillator 
\cite{hinterding,KS}.

Including negative eigenvalues of $x$, i.e. taking the direct sum of 
two irreducible representations of the algebra generated by $X$, $P$ and 
$U$, is also a possibility to obtain a Hilbert space, on which $X$ and $P$ 
are represented by self-adjoint operators \cite{FLW}

We now turn to the unbounded part of the spectrum, i.e.
$\ve=q^{2\ga-2m}, m\in\mathbb{Z}$. To connect our solutions with
the results of \cite{Koorn} (Eqn (4.2)), it is necessary, to set
$\xi^2\la^2=q$ or, equivalently:
\begin{equation}
\label{c-ga-ident}
  c=q^{-2\ga}.
\end{equation}
For the functions ${_1}\vp_1(-cq^{-4\nu+4p};q^{-4\nu-4};q^{-4},c^{-1}xq^{-4p-4})$,  
with $\nu=\pm\ha$, the square of the norm is according to \cite{Koorn}:
\begin{equation}
  \delta_{pr}cq^{4p+2}\frac{(-c^{-1}q^{-4p-4},-c^{-1}q^{-2};q^{-4})_{\infty}}{(-c^{-1}q^{-4p-6},-cq^{-2},-c,-c^{-1}q^{-4};q^{-4})_{\infty}}
  \left(\frac{(q^{-4},-cq^{-2};q^{-4})_{\infty}}{(q^{-2};q^{-4})_{\infty}}\right)^2
\end{equation}
and
\begin{equation}
  \delta_{pr}cq^{4p-2}\frac{(-c^{-1}q^{-4p-4},-c^{-1}q^{2};q^{-4})_{\infty}}{(-c^{-1}q^{-4p-2},-cq^{-6},-c,-c^{-1}q^{-4};q^{-4})_{\infty}}
  \left(\frac{(q^{-4},-cq^{-6};q^{-4})_{\infty}}{(q^{-6};q^{-4})_{\infty}}\right)^2
\end{equation}
If we take (\ref{c-ga-ident}) into account we find for the
constants $C_e, C_o$ in (\ref{gesLsg}): In the case
$\ve=q^{2\ga-4p-2}=c^{-1}q^{-4p-2}$
\begin{equation}
  C_e=(-c)^pq^{4p^2+2p}\frac{(q^{-2},q^{-4})_{\infty}}{(-c^{-1}q^{-4p-4};q^{-4})_{\infty}},\qquad
  C_o=0
\end{equation}
and for  $\ve=q^{2\ga-4p}=c^{-1}q^{-4p}$
\begin{equation}
  C_e=0,\qquad
  C_o=-(-c)^pq^{4p^2-2p}\frac{(q^{-6},q^{-4})_{\infty}}{(-c^{-1}q^{-4p-4};q^{-4})_{\infty}}
\end{equation}

Now, take $m=2p+1$ in the first case and $m=2p$ in the second case. Such that 
$\ve=q^{2\ga-2m}$ with $m\in\mathbb{Z}$.  
Putting together these results, the Norm of the combined solution (\ref{gesLsg}) becomes 
for both cases
\begin{equation}
  \delta_{pr}c^{m}q^{2m^2}\frac{M_c}{(-c^{-1}q^{-2m-2};q^{-2})_{\infty}}, 
\end{equation}
with $M_c=\frac{(q^{-4},q^{-4},-c^{-1}q^{-2},-cq^{-2};q^{_4})_{\infty}}{(-c,-c^{-1}q^{-4};q^{-4})_{\infty}}$.

We define
\begin{equation}
  \tilde{k}_m(x)=(-x)^{m-\ga}\sqrt{c}^{\ga-m}{_2}\vp_1\left(-q^{2m-2\ga-2},-q^{2m-2\ga};0;q^{-4},-\frac{q^{-4}}{x^2}\right);
\end{equation}
if we extend the lattice to negative values as above, 
these functions are orthogonal:
\begin{equation}
  \sum_{k=-\infty,\si=\pm}^{+\infty}\frac{\tilde{k}_m(\si\sqrt{c}q^{-2k})\tilde{k}_n(\si\sqrt{c}q^{-2k})}{(-cq^{-4k};q^{-4})_{\infty}}q^{-2k}=2M_c\frac{c^mq^{2m^2}}{(-c^{-1}q^{-2m-2},q^{-2})_{\infty}}\delta_{mn},
\end{equation}
Also the scalar product of two functions $\tilde{h}_n$ and $\tilde{k}_m$ 
vanishes. The results of \cite{Koorn} imply that
the set of functions $\tilde{h}_n$, $n\in\mathbb{N}$ together 
with the functions $\tilde{k}_m$, $m\in\mathbb{Z}$ form a basis of the 
Hilbert space.

The $q$-Heisenberg algebra (\ref{XPURel}) is represented on the space of 
square integrable functions on the set
\begin{equation}
  \Lambda=\{\sigma\xi_oq^n| \sigma=\pm 1, n\in\mathbb{Z}\},\; \xi_o\in\mathbb{R}.
\end{equation}
With the scalar product given by a Jackson integral:
\begin{equation}
  (f,g)=\xi_o(q-q^{-1})\sum_{\sigma,n}\overline{f(\sigma\xi_oq^n)}g(\sigma\xi_oq^n)q^n.
\end{equation}
The operators $X$, $P$, $U$ act according to (\ref{XPU-realization}). 
$X$ and $P$ are essentially self-adjoint, $U$ is unitary. A representation 
is characterized by $\xi_o\in [1,q[$.

The $q$-difference equation (\ref{f-EV-eqn}) for the eigenvalues of the 
Hamilton operator $H=\ac a$ only connects even and odd lattice points among 
them self. That means, there is a twofold degeneracy in the spectrum. 

Depending on the lattice that we consider (even or odd), the 
parameter $c$ has to be related to $\xi_o$ in different ways. 
For even lattice points, $\sigma \xi_o q^{2n}$, we have $\sqrt{c}=\xi_o$. 
For odd lattice points, $\sigma \xi_o q^{2n+1}$, we find $q\sqrt{c}=\xi_o$ 
and a shift in $\ga$ occurs, which however does not change the spectrum. 

\section{Summary}

It is well known, that in general spectra of Hamilton operators related to deformed 
oscillator algebras consist of several different parts \cite{AOd}. The example 
discussed in this paper shows, that in order to get self-adjoint representations, 
one has to take into account all parts. 
It would be very interesting to see, whether these parts play a role in a quantum 
field theory, that is constructed with the aid of a deformed oscillator algebra. 
 
In our specific case it turned out, that a nice way to get a self-adjoint representation 
is to consider not only positve lattice points, which already form a irreducible 
representation of the algebra (\ref{XPURel}) , but also negative. 
This has also been done in \cite{FLW} in order to get a representation, such that 
$X$ and $P$ are both essentially self-adjoint. 
It also was shown, that in the unbounded part $\frac{1+q^{2\ga-2m}}{1-q^{-2}},m\in\mathbb{Z},\; \ga\in\mathbb{R}$ of the spectrum one $\ga$ is singled out by the chosen lattice, cf. 
(\ref{c-ga-ident}).

As explained in \cite{Koorn}, these results also show, that the momentum problem 
associated to the weight function $\frac{1}{(-x^2;q^{-4})_{\infty}}$, which appears 
in the orthogonality relations, e.g. Eq. (\ref{Hermite-ortho}), is indetermined. 
Each of the functions $\tilde{k}_s$ is bounded: $|\tilde{k}_s(\pm\sqrt{c}q^{-2n})|<C$, 
for all $n$. Since the functions $\tilde{k}_s$ are orthogonal to the Hermite polynomials 
$\tilde{h}_m$ and therefore to all polynomials, the moments will not change if one 
uses for example $\frac{1+C^{-1}\tilde{k}_s(x)}{(-x^2;q^{-4})_{\infty}}$ as weight function. 
From another point of view, see e.g. \cite{Simon}, this happens, because the operator $X$ 
is not self-adjoint in the space spanned by the $q$-Hermite polynomials together 
with the scalar product, that is given by the moments. In some sense one may 
interpret specifying $\ga$ as choosing a self-adjoint extension.


\end{document}